\documentclass[12pt, reqno]{amsart}
\usepackage{amsmath}
\usepackage{amsthm}
\usepackage{amssymb}
\usepackage{enumerate}
\usepackage{mathrsfs}
\usepackage{graphicx}
\usepackage{graphicx}
\usepackage{amssymb,amsfonts}
\usepackage{amsmath}
\usepackage[colorlinks,linkcolor=blue,anchorcolor=blue,citecolor=blue]{hyperref}
\usepackage[numbers,sort&compress]{natbib}
\usepackage{marvosym}
\usepackage{epstopdf}
\usepackage{caption}
\usepackage{multicol}
\usepackage{multirow}
\usepackage{makebox}
\usepackage{subfigure}
\usepackage{dsfont}
\numberwithin{equation}{section}
\newtheorem{theo}{Theorem} 

\newtheorem{lem}{Lemma}

\begin{document}

\title {$\textrm{SL}(n)$ covariant matrix-valued valuations on Orlicz spaces}

\author[Zeng and Lan]{Chunna Zeng$^{*}$ and Yu Lan}

\address{1.School of Mathematical Sciences, Chongqing Normal University, Chongqing 401331, People's Republic of China; Institut f\"{u}r Diskrete Mathematik und Geometrie, Technische Universit\"{a}t Wien, Wiedner Hauptstra\ss e 8--10/1046, 1040 Wien, Austria }
\email{zengchn@163.com}

\address{2.School of Mathematics and Sciences, Chongqing Normal University, Chongqing 401332, People's Republic of China}
\email{18181567221@163.com}

\thanks{The first author is supported by the Major Special Project of NSFC (Grant No. 12141101), the Young Top-Talent program of Chongqing (Grant No. CQYC2021059145), the Science and Technology Research Program of Chongqing Municipal Education Commission (Grant No. KJZD-K202200509) and Natural Science Foundation Project of Chongqing (Grant No. CSTB2024NSCQ-MSX0937)}
\thanks{{*}Corresponding author: Chunna Zeng}
\begin{abstract}
All continuous,  $\mathrm{SL}(n)$ covariant valuations on Orlicz spaces are completely classified without any symmetric assumptions.  It is shown that the moment matrix is the only such valuation if $n\geq3$, while a new functional shows up in dimension two.
\end{abstract}

\maketitle

\section{Introduction}
Let $\mathcal{S}$ be the subset of a lattice $(\mathcal{L},\cup,\cap)$ and let $\mathcal{A}$ be an abelian semigroup.
 An operator $\mathrm{Y}:\mathcal{S}\rightarrow\langle \mathcal{A},+\rangle$ is called a valuation if
\begin{equation*}
\mathrm{Y}(Q_{1})+\mathrm{Y}(Q_{2})=\mathrm{Y}(Q_{1}\cup Q_{2})+\mathrm{Y}(Q_{1}\cap Q_{2})
\end{equation*}
for all $Q_{1},Q_{2} \in \mathcal{S}$ and $Q_{1}\cup Q_{2}, Q_{1}\cap Q_{2} \in \mathcal{S}$. Here, the operations $\cup$ and $\cap$ are union and intersection, respectively. Since Dehn's solution of Hilbert's Third Problem in 1901, many excellent results have emerged in the study of valuations. One of the most well-known results is Hadwiger's characterization theorem, which classified all continuous and rigid motion invariant valuations on the space of convex bodies (compact convex sets). Hadwiger's theorem has wide applications in integral geometry and geometric probability. Ludwig has made significant and profound contributions to the field of valuation theory. For a review of body-valued valuations, their historical development and recent contributions on convex bodies, star-shaped sets, manifolds, and related topics, see \cite{Alesker1, Alesker2, Alesker3, Alesker4, Klain1, Klain2, Klain3, Klain-Rota, L-Y-Z1, Ludwig1, Ludwig2, Ludwig3, Ludwig4, Ludwig8, Ludwig9, McMullen1, McMullen2, Parapatits1, Parapatits2, Wannerer, Wang1, Zeng-Ma, Zeng-Zhou}.
\par The results related to body-valued valuations have facilitated further exploration of valuations on function spaces. Suppose that $\mathcal{S}$ is a space of real-valued functions, and $\mathbb{M}^n$ denotes the set of $n \times n$-matrices. 
An operator $\Psi:\mathcal{S}\rightarrow \langle\mathbb{M}^n,+\rangle$ is called a valuation if
\begin{equation*}
\Psi(h_{1})+\Psi(h_{2})=\Psi(h_{1}\vee h_{2})+\Psi(h_{1}\wedge h_{2})
\end{equation*}
for all $h_{1}, h_{2}\in\mathcal{S}$ and $h_{1}\vee h_{2}, h_{1}\wedge h_{2} \in \mathcal{S}$. Here $h_{1}\vee h_{2}$ and $h_{1}\wedge h_{2}$ represent the pointwise maximum and pointwise minimum of $h_{1}$ and $h_{2}$, respectively. Recently, the study of valuations on function spaces has seen substantial progress, including $L^{p}$-spaces \cite{Ludwig6, Tsang1, Tsang2, Wang2}, Sobolev spaces \cite{Ludwig5, Ludwig7}, the space of Lipschitz functions \cite{Colesanti1, Colesanti2}, the spaces of functions of bounded variation \cite{Wang}, Orlicz spaces \cite{Kone} and the space of convex functions \cite{Colesanti3, Colesanti4, Ludwig10}, among others. In \cite{Tsang1}, a complete classification of continuous and translation invariant real-valued valuations on $L^{p}(\mathbb{R}^n)$ was established by Tsang. Based on Tsang's work, Ludwig \cite{Ludwig6} classified the continuous, $\mathrm{SL}(n)$ covariant symmetric matrix-valued valuations on $L^1(\mathbb{R}^n,|x|^2dx) $. Furthermore, Wang proved that the symmetry assumption in \cite{Ludwig6} can be removed. Kone \cite{Kone} generalized Tsang’s results from  $L^{p}$-spaces to Orlicz spaces, and also demonstrated that Klain’s results \cite{Klain2} can be extended to $L^{\phi}$-star valuations.

Let $\phi$ be a Young function. Denote the Orlicz space by $L^{\phi}(\mathbb{R}^n,|x|^2dx)$, which is the space of measurable
functions $h:\mathbb{R}^n \rightarrow \mathbb{R}$ such that $\int_{\mathbb{R}^n}\phi(\delta|h(x)|)|x|^2dx<\infty$ for any $\delta>0$, where $|x|$ denotes the Euclidean norm of $x$ in $\mathbb{R}^n$. An operator $\Psi: L^{\phi}(\mathbb{R}^n,|x|^2dx) \rightarrow \langle \mathbb{M}^n, +\rangle$ is said to be $\mathrm{SL}(n)$ covariant if
\begin{equation*}
\Psi(h \circ \vartheta^{-1})=\vartheta\Psi(h)\vartheta^{t}
\end{equation*}
for $\vartheta \in \mathrm{SL}(n)$.  Define the moment matrix $\mathrm{K}(h)$ of a measurable function $h:\mathbb{R}^n \rightarrow \mathbb{R}$, which is the $n \times n$-matrix with (not necessarily finite) entries given by
\begin{equation*}
\mathrm{K}_{ij}(h)=\int_{\mathbb{R}^n}h(x)x_{i}x_{j}dx.
\end{equation*}

 The main aim of this paper is to establish a complete classification of continuous, $\mathrm{SL}(n)$ covariant matrix-valued valuations on Orlicz spaces. Let $C_{0}(\mathbb{R})$ stand for the collection of continuous functions $\xi:\mathbb{R}\rightarrow\mathbb{R}$, where $\xi(0)=0.$
 Let $C_{\phi}(\mathbb{R})$ represent the collection of continuous functions
$\xi:\mathbb{R}\rightarrow\mathbb{R}$ such that $|\xi(\beta)|\leq\gamma\phi(|\beta|)$ for some $\gamma\geq0$ and all $\beta\in\mathbb{R},$  where $\phi$ is a Young function. The detailed definition of the $\Delta_2$-condition will be introduced in Section 2.
Our main results are the matrix-valued analogs of Kone’s characterization of real-valued valuations (see \cite{Kone}).

\begin{theo}\label{theo2}
Let $\phi$ be a Young function that satisfies the $\Delta_2$-condition.
An operator $\Psi:L^{\phi}(\mathbb{R}^2,|x|^2dx)\rightarrow \langle\mathbb{M}^2, +\rangle$
is a continuous and $\mathrm{SL}(2)$ covariant valuation if and only if there exists a functional $\xi\in C_{\phi}(\mathbb{R})$ and constant $c\in\mathbb{R}$ such that
\begin{equation*}
\Psi(h)=\mathrm{K}(\xi\circ h)+c\rho_{\pi/2}
\end{equation*}
for $h\in L^{\phi}(\mathbb{R}^2,|x|^2dx)$. Here $\rho_{\pi/2}$ denotes the counter-clockwise rotation  by $\pi/2$ in $\mathbb{R}^2$.
\end{theo}

\begin{theo}\label{theo3}
Let $n\geq 3$ and $\phi$ be a Young function that satisfies the $\Delta_2$-condition. An operator $\Psi:L^{\phi}(\mathbb{R}^n,|x|^2dx)\rightarrow \langle\mathbb{M}^n, +\rangle$
is a continuous and $\mathrm{SL}(n)$ covariant valuation if and only if there exists a functional $\xi\in C_{\phi}(\mathbb{R})$ such that
\begin{equation*}
\Psi(h)=\mathrm{K}(\xi\circ h)
\end{equation*}
for $h\in L^{\phi}(\mathbb{R}^n,|x|^2dx)$.
\end{theo}
\section{~Preliminaries}
\subsection{Measure Spaces}

\par A measure space is a triple $(X,\mathcal{F},\mu)$, where
\begin{enumerate}
\item $X$ is a set referred to as the sample space;
\item $\mathcal{F}$ is a $\sigma$-algebra on $X$;
\item $\mu:\mathcal{F}\rightarrow[0,\infty)$ is a measure defined on $\mathcal{F}$.
\end{enumerate}
\par In this paper we give two measurable spaces: $(\mathbb{R}^n,\mathscr{B},\lambda)$ and $(\mathbb{R}^n,\mathscr{B}, |x|^2dx)$, where $\mathscr{B}$ represents the collection of Lebesgue measurable sets in $\mathbb{R}^n$, $\lambda$ denotes the standard Lebesgue measure and $|x|^2dx$ is weighted Lebesgue measure in $\mathbb{R}^n$. For convenience, we denote $\mu_{n}(M)=\int_{M}|x|^2dx$, where $M \subset \mathbb{R}^n$.

\subsection{Background material on Orlicz spaces.}
Young functions play a crucial role in Orlicz spaces.
A Young function $\phi$ is a non-negative, real-valued function  defined on $[0,\infty)$, which can be expressed as 
\begin{equation*}
\phi(t)=\int_{0}^{t}\varphi(s)ds \,\,\mathrm{for} \,\ t\geq0.
\end{equation*}
 The function $\varphi: [0, \infty) \rightarrow [0, \infty)$ is non-decreasing, right-continuous, with $\varphi(0) = 0$ and $\lim_{s \rightarrow \infty} \varphi(s)=\infty$.
Young function $\phi$ is convex, non-negative, continuous, and increasing on $[0,\infty).$ Moreover, $\phi(0)=0$ and $\lim_{t\rightarrow\infty}\phi(t)=\infty$. The function $\phi^{*}$ is the complementary function to $\phi$, defined by
\begin{equation*}
\phi^{*}(t)=\int_{0}^{t}\psi(s)ds,
\end{equation*}
where $\psi(t)=\mathrm{sup}\{s:\varphi(s)\leq t\}$
for $t\geq0$. The pair $(\phi,\phi^{*})$ is referred to as complementary Young functions. Additionally, the complementary function $\phi^{*}$  of $\phi$ is also a Young function. The following lemma provides further properties of the Young function $\phi$.
\begin{lem}\label{lem1}\cite{Kone}
If $\phi$ is a Young function, then
\begin{equation*}
\lim\limits_{t\rightarrow\infty}\frac{\phi(t)}{t}=\infty \,\,and \,\,\lim\limits_{t\rightarrow\infty}\frac{\phi^{-1}(t)}{t}=0
\end{equation*}
for all $t\geq0$.
\end{lem}

 \par The Orlicz class $ L_{\dagger}^{\phi}(\mathbb{R}^n, |x|^2dx)$ is defined as
\begin{equation*}
L_{\dagger}^{\phi}(\mathbb{R}^n, |x|^2dx)=\{h:\mathbb{R}^n\rightarrow\mathbb{R},\,\mathrm{measurable}, \,\rho_{\phi}(h)<\infty\},
\end{equation*}
where \begin{equation*}
\rho_{\phi}(h)=\int_{\mathbb{R}^n}\phi(|h(x)|)|x|^2 dx.
\end{equation*}
And the Orlicz space $L^{\phi}(\mathbb{R}^n, |x|^2dx)$ is defined as
\begin{equation*}
L^{\phi}(\mathbb{R}^n, |x|^2dx)=\{h:\mathbb{R}^n\rightarrow\mathbb{R},\,\mathrm{measurable}, \,\rho_{\phi}(\delta h)<\infty\, \,\mathrm{for}\,\, \mathrm{any} \,\, \delta>0\}. 
\end{equation*}
A Young function $\phi$  satisfies the $\Delta_{2}$-condition
if there exist $T\geq0$ and a positive constant $k$ such that  
\begin{equation*}
\phi(2t)\leq k\phi(t)
\end{equation*}
for all $t\geq T$.
If $\phi$ satisfies $\Delta_{2}$-condition, then the Orlicz class is equivalent to the Orlicz space.
\par 
The Orlicz norm $\| h\|_{\phi}$ of $h \in L^{\phi}(\mathbb{R}^n, |x|^2dx)$  is given by
\begin{equation*}
\| h\|_{\phi}=\mathrm{sup}\left\{\int_{\mathbb{R}^n}|h(x)g(x)||x|^2 dx: g\in L^{\phi^{*}}_{\dagger}(\mathbb{R}^n, |x|^2dx), \rho_{\phi^*}(g)\leq1\right\}.
\end{equation*}
The functional $\|\cdot\|_{\phi}: L^{\phi}(\mathbb{R}^n, |x|^2dx)\rightarrow\mathbb{R}$ is a semi-norm. Since we identify $h$ and $g$ when $\| h-g\|_{\phi}=0$, this turns $\| \cdot\|_{\phi}$ into a norm, making $L^{\phi}(\mathbb{R}^n, |x|^2dx)$  a normed linear space. 
For a sequence of functions $\{h_k\}$ in $L^{\phi}(\mathbb{R}^n, |x|^2dx)$,   $h_k$ converges to  $h$ in $L^{\phi}(\mathbb{R}^n, |x|^2dx)$ if 
 \begin{equation*}
 \lim\limits_{k\rightarrow \infty} \|h_k - h\|_{\phi} = 0.    
 \end{equation*}
The following result involves the convergence of a sequence of functions in Orlicz spaces.
 \begin{lem}\label{lem6}\cite{Kone}
  Let $\{h_{i}\}\subset L^{\phi}(\mathbb{R}^n, |x|^2dx)$ be a sequence and $h\in L^{\phi}(\mathbb{R}^n, |x|^2dx)$ such that $\Vert h_{i}-h\Vert_{\phi}\rightarrow0$. Then, there exist a subsequence $\{h_{k}\}$ of $\{h_{i}\}$ and a functional $g^*\in L^{\phi}(\mathbb{R}^n, |x|^2dx)$ such that
 
 $(1)\ h_{k}\rightarrow h$ a.e.; 
 
$(2)\ |h_{k}|\leq g^*$ a.e.  for all $k$.
\end{lem}

Convergence in norm provides a foundation for analyzing continuity in normed linear spaces.
An operator $\Psi: L^{\phi}(\mathbb{R}^n, |x|^2dx) \rightarrow \langle\mathbb{M}^n,+ \rangle$ is said to be continuous if, for every sequence $\{h_k\}$ converging to $h$ in $L^{\phi}(\mathbb{R}^n, |x|^2dx)$, we have $\Psi(h_k) \rightarrow \Psi(h)$.

\par Next we introduce the subspace $ L^{\phi}_{a}(\mathbb{R}^n, |x|^2dx)$ of $L^{\phi}(\mathbb{R}^n, |x|^2dx) $. 
The annulus $A[r, R)$ represents the region between two concentric spheres and is given by
\begin{equation*}
A[r,R)=\{x\in\mathbb{R}^n: r\leq |x|<R\},
\end{equation*}
where $r$ and $R$ are real numbers satisfying $0\leq r<R,$ and $|x|$ 
 is the distance from the origin to  $x$.
The subspace $ L^{\phi}_{a}(\mathbb{R}^n, |x|^2dx) $ is defined by
\begin{equation*}
L^{\phi}_{a}(\mathbb{R}^n, |x|^2dx) = \{ h \in L^{\phi}(\mathbb{R}^n, |x|^2dx) : \text{supp}(h) \subseteq A[r, R) \},
\end{equation*}
where the support of $h,$ denoted by supp$(h)$, is the set of points in $\mathbb{R}^n$ such that $h\neq0$.
The definition implies that $ h $ vanishes outside the annulus $ A[r, R) $.
\par The indicator functions and simple functions play crucial roles in our subsequent proofs. Define $\chi_{M}$ as the indicator function of $M \subset \mathbb{R}^n$, that is, $\chi_{M}(x)=1$ for $x\in M$ and $\chi_{M}(x)=0$ otherwise. If $s$ is a simple function, then there exist pairwise disjoint sets $M_{i} \in\mathscr{B}$ that have finite measure and scalars $\alpha_{i}$,  $i=1,\ldots,l$ such that
\begin{equation*}
s=\sum\limits_{i=1}\limits^{l}\alpha_{i}\chi_{M_{i}}(x).
\end{equation*}
Using the definition of indicator functions and the Orlicz norm, Kone obtained the following result.
\begin{lem}\label{lem3}\cite{Kone}
Let $M\subset\mathbb{R}^n$. If $0<\mu_{n}(M)<\infty$, then
\begin{equation*}
\|\chi_{M}\|_{\phi}=\phi^{*^{-1}}\left(\frac{1}{\mu_{n}(M)}\right)\mu_{n}(M).
\end{equation*}
 \end{lem}
 
\subsection{Background material on convex polytopes}

\par Let $e_{1},\ldots,e_{n}$ denote the standard basis vectors in $\mathbb{R}^n$.
Let $[u_{1},\ldots,u_{k}]$ represent the convex hull of points $u_{1},\ldots,u_{k}\in\mathbb{R}^n$. Denote $\mathcal{P}^n$ by the collection of compact convex polytopes, endowed with the topology induced by the Hausdorff metric.
For $Q\in\mathcal{P}^{n}$, the moment matrix $\mathrm{M}(Q)$ is defined as  
 \begin{equation*}
\mathrm{M}(Q)=\int_{Q}x_{i}x_{j}dx.
 \end{equation*}
An operator $\mathrm{Y}:\mathcal{P}^n\rightarrow \langle\mathbb{M}^{n}, +\rangle$ is $\mathrm{SL}(n)$ covariant if 
\begin{equation*}
\mathrm{Y}(\vartheta Q)=\vartheta \mathrm{Y}(Q)\vartheta^t
\end{equation*}
for any $Q\in \mathcal{P}^n$ and $\vartheta\in \mathrm{SL}(n)$. An operator $\mathrm{Y}:\mathcal{P}^n\rightarrow\langle \mathbb{M}^n,+\rangle$ is called weak simple if there is a constant matrix $\mathrm{M}_{0}\in\mathbb{M}^n$ satisfying $\mathrm{Y}(Q)=\mathrm{M}_{0}$, where $Q\in\mathcal{P}^n$ with dim $Q\leq n-1$. If $\mathrm{M}_{0}=\textbf{0}$, where $\textbf{0}$ denotes the zero matrix, then the operator $\mathrm{Y}$ is simple.
Note that the moment matrix $\mathrm{M}$ is a simple, $\mathrm{SL}(n)$ covariant, and continuous function.

\par 
A closed cube $\mathcal{C}$ in $\mathbb{R}^n$  is characterized by the set 
\begin{equation*}
\mathcal{C}=\{(x_{1},\ldots,x_{n}): c_{i}\leq x_{i}\leq d_{i}, x_{i}, c_{i}, d_{i}\in\mathbb{R}\},
\end{equation*}
where $\vert c_{i}-d_{i}\vert= \vert c_{j}-d_{j}\vert$,  $1\leq i, j\leq n$. For two sets $M_1, M_2\subseteq\mathbb{R}^n$, their symmetric difference is expressed as  
 \begin{equation*}
 M_1\bigtriangleup M_2=M_1\backslash M_2 \cup M_2\backslash M_1.
 \end{equation*}

\begin{lem}\label{lem2}\cite{Tsang1}
Given that $\varepsilon>0$ and $M\in\mathscr{B}$ with $\lambda(M)<\infty$, there is a sequence of finite closed cubes $C_{1}, C_{2},\ldots, C_{l}\subseteq\mathbb{R}^n$ such that $\lambda(M\bigtriangleup \bigcup\limits_{i=1}\limits^{l} C_{i})<\varepsilon$.    
\end{lem}


\section{Proof of main theorems}
\begin{theo}\label{theo2} Let $n\geq2$. If $\xi \in C_{0}(\mathbb{R})$, and for  any $\delta>0$, there exists a real number $\gamma_{\delta}\geq 0$ such that $|\xi(\beta)|\leq\gamma_{\delta}\phi(\delta|\beta|)$ for $\beta \in \mathbb{R}$, then the operator $\Psi:L^{\phi}(\mathbb{R}^n, |x|^2dx)\rightarrow \langle\mathbb{M}^n,+\rangle$ defined by
\begin{equation*}
\Psi(h)=\mathrm{K}(\xi\circ h),
\end{equation*}
is a continuous and $\mathrm{SL}(n)$ covariant valuation on $L^{\phi}(\mathbb{R}^n, |x|^2dx)$.
\end{theo}
\emph{Proof.}
Since $| \xi(\beta)|\leq\gamma_{\delta}\phi(\delta|\beta|)$ for $\beta\in\mathbb{R}$, then 
\begin{align*}
|\mathrm{K}_{ij}(\xi\circ h)|&=\left|\int_{\mathbb{R}^{n}}(\xi\circ h)(x)x_{i}x_{j}dx\right|\leq \int_{\mathbb{R}^{n}}|(\xi\circ h)(x)||x|^2dx\\
&\leq \gamma_{\delta}\int_{\mathbb{R}^{n}}\phi(\delta|h(x)|)|x|^2dx<\infty.
\end{align*}
It shows that $|\Psi(h)|<\infty$ for $h\in L^{\phi}(\mathbb{R}^n, |x|^2dx)$. Furthermore, from $\xi(0)=0$, it follows that $\Psi(0)=\mathrm{\textbf{0}}.$
Since \begin{equation*}
\mathrm{K}( h)=\int_{\mathbb{R}^n}xx^t| h(x)|dx
\end{equation*} 
for any $h_{1},h_{2}\in L^{\phi}(\mathbb{R}^n, |x|^2dx)$, we have
\begin{align*}
\Psi(h_{1}\vee h_{2})+\Psi(h_{1}\wedge h_{2})&=\int_{\mathbb{R}^n}xx^{t}\vert\xi\circ(h_{1}\vee h_{2})\vert dx+\int_{\mathbb{R}^n}xx^{t}\vert\xi\circ(h_{1}\wedge h_{2})\vert dx\\
&=\int_{\{h_{1}\geq h_{2}\}}xx^{t}\vert\xi\circ h_{1}\vert dx+\int_{\{h_{1}<h_{2}\}}xx^{t}\vert\xi\circ h_{2}\vert dx\\
&\, \, \, \, \, +\int_{\{h_{1}\geq h_{2}\}}xx^{t}\vert\xi\circ h_{2}\vert dx+\int_{\{h_{1}<h_{2}\}}xx^{t}\vert\xi\circ h_{1}\vert dx\\
&=\int_{\mathbb{R}^n}xx^{t}\vert\xi\circ h_{1}\vert dx+\int_{\mathbb{R}^n}xx^{t}\vert\xi\circ h_{2}\vert dx\\
&=\Psi(h_{1})+\Psi(h_{2}),
\end{align*}
which implies that $\Psi$ is a valuation on $L^{\phi}(\mathbb{R}^n, |x|^2dx)$.
\par 
Due to the $\mathrm{SL}(n)$  covariance of $\mathrm{K}$, then
\begin{equation*}
\Psi(h\circ\vartheta^{-1})=\mathrm{K}(\xi\circ(h\circ\vartheta^{-1}))=\mathrm{K}((\xi\circ h)\circ\vartheta^{-1})=\vartheta \mathrm{K}(\xi\circ h)\vartheta^{t}=\vartheta \Psi(h)\vartheta^{t}
\end{equation*}
for all $\vartheta\in \mathrm{SL}(n)$. Thus, $\Psi$ is $\mathrm{SL}(n)$ covariant.

\par
 To establish the continuity of $\Psi$, we need to prove that $\Psi(h_{i})\rightarrow\Psi(h)$ as $h_{i}\rightarrow h$ in $L^{\phi}(\mathbb{R}^n, |x|^2dx)$.
Suppose that $\{h_{i}\}$ is a sequence of functions in $L^{\phi}(\mathbb{R}^n, |x|^2dx),$ and $h_{i}\rightarrow h$ in $L^{\phi}(\mathbb{R}^n, |x|^2dx)$. By Lemma \ref{lem6}, there is a subsequence $\{h_{k}\}$ of $\{h_{i}\}$ satisfying $h_{k}\rightarrow h$ a.e., and $g^*\in L^{\phi}(\mathbb{R}^n, |x|^2dx)$ such that $|h_{k}|\leq g^*$ a.e. for all $k$. Since $\xi$ is continuous, it follows that $(\xi\circ h_{k})(x)x_{i}x_{j}\rightarrow(\xi\circ h) (x)x_{i}x_{j}$ a.e. as $k\rightarrow\infty$. Additionally, 
since  $|\xi(\beta)|\leq\gamma_{\delta}\phi(\delta\vert \beta\vert)$ for $\beta\in \mathbb R$, we have
\begin{equation*}
|(\xi\circ h_{k})(x)x_{i}x_{j}|\leq|(\xi\circ h_{k})(x)||x|^2
\leq\gamma_{\delta}\phi(\delta|h_{k}(x)|)|x|^2\leq\gamma_{\delta}\phi(\delta g^*)|x|^2,\, a.e.
\end{equation*}
Observe that $\gamma_{\delta}\phi(\delta g^*)|x|^2$ is integrable, by the dominated convergence theorem, we obtain
\begin{equation*}
\lim\limits_{k\rightarrow\infty}\int_{\mathbb{R}^n}(\xi\circ h_{k})(x)x_{i}x_{j} dx=\int_{\mathbb{R}^n}(\xi\circ h)(x)x_{i}x_{j}dx.
\end{equation*}
Thus, $\Psi(h_{k})\rightarrow \Psi(h)$ as $h_{k}\rightarrow h$ in $L^{\phi}(\mathbb{R}^n, |x|^2dx)$, which demonstrates that $\Psi$ is continuous.  \qed


\begin{lem}\label{lem3}
Suppose that $\Psi:L^{\phi}(\mathbb{R}^2, |x|^2dx)\rightarrow \langle \mathbb{M}^2, + \rangle$ is an $\mathrm{SL}(2)$ covariant. Then there is a real constant $c$ such that
\begin{equation*}
\Psi(0)=c\rho_{\pi/2}.
\end{equation*}
 Here, $\rho_{\pi/2}$ denotes the counter-clockwise rotation  by $\pi/2$ in $\mathbb{R}^2$.
\end{lem}

\par \emph{Proof.} 
Let $\vartheta_{1}, \vartheta_{2}\in \textrm{SL}(2)$, where $\vartheta_{1}$=diag$(k,k^{-1})$ with $k\in\mathbb R \backslash \{0\}$
and $\vartheta_{2}=$
$\begin{pmatrix}
1 & 0\\
1 & 1
\end{pmatrix}$.
The  $\textrm{SL}(2)$ covariance of $\Psi$ implies that
\begin{equation}\label{6}
 \Psi(0)= \Psi(0\circ\vartheta_{i}^{-1})=\vartheta_{i}\Psi(0)\vartheta_{i}^{t}
\end{equation}
 for $i=1,2$.
Given that $\Psi(0)=(a_{ij})_{2\times2}$,
by substituting $\vartheta_1$ into (\ref{6}), we obtain
 \begin{equation*}
\begin{pmatrix}
 a_{11} & a_{12}\\
a_{21} & a_{22}
 \end{pmatrix}=
 \begin{pmatrix}
 k^{2}a_{11} & a_{12}\\
 a_{21} & k^{-2}a_{22}
 \end{pmatrix}
 \end{equation*}
 for any $k\in\mathbb R \backslash \{0\}$. Hence, $a_{11}=0$ and $a_{22}=0$.
 Taking $\vartheta_2$ in (\ref{6}) we have
\begin{equation*}
\begin{pmatrix}
 0 & a_{12}\\
a_{21} & 0
 \end{pmatrix}=
 \begin{pmatrix}
 0 & a_{12}\\
 a_{21} & a_{12}+a_{21} \end{pmatrix},
 \end{equation*}
 so $a_{12}=-a_{21}$.
 \qed

\begin{lem}\label{lem4}
Suppose that $n\geq3$. An operator $\Psi:L^{\phi}(\mathbb{R}^n, |x|^2dx)\rightarrow\langle \mathbb{M}^n ,+\rangle$ is $\mathrm{SL}(n)$ covariant, then 
\begin{equation*}
\Psi(0)=\mathrm{\textbf{0}}.
\end{equation*}
\end{lem}

\emph{Proof.} 
For a real number $k\in\mathbb{R}\backslash\{0\}$, let
$\vartheta_{n}=\mathrm{diag}(k^{n-1},k^{-1},\ldots,k^{-1}),$
then $\vartheta_{n}\in\mathrm{SL}(n)$. Assume that
$\Psi(0)=\left(a_{ij}\right)_{n\times n}$. Because $\Psi$ has the property of  $\mathrm{SL}(n)$ covariant, then $\Psi(0)= \Psi(0\circ\vartheta_{n}^{-1})=\vartheta_{n}\Psi(0)\vartheta_{n}^{t}.$ A direct calculation gives
 \begin{equation}\label{003}
 a_{ij}=k^{-2}a_{ij},\ \ a_{i1}=k^{n-2}a_{i1}, \ \  a_{1i}=k^{n-2}a_{1i},\ \ \ 2 \leq i, j\leq n;
\end{equation} 
and
 \begin{equation}\label{004}
a_{11}=k^{2n-2}a_{11}.
\end{equation} 
By combining (\ref{003}) and (\ref{004}), and the arbitrariness of 
$k$, we obtain
 \begin{equation}\label{005}
  a_{ij}=0,\ \ a_{i1}=a_{1i}=0,\ \ a_{11}=0,\ \ 2 \leq i, j\leq n.
\end{equation}
Thus, $\Psi(0)=\textbf{0}.$
\qed

Recently, Wang \cite{Wang2} established that the moment matrix is the only classification of continuous and $\mathrm{SL}(n)$ covariant matrix-valued valuations on convex polytopes when $n\geq 3$, while $\rho_{\pi/2}$ appears in the $2$-dimensional case.

\begin{lem}\label{lem8}\cite{Wang2}
An operator $\mathrm{Y}:\mathcal{P}^2\rightarrow \langle \mathbb{M}^2 ,+\rangle$ is a weak simple, continuous and $\mathrm{SL}(2)$ covariant valuation if and only if there are real constants $c_{1},c_{2}$ such that
\begin{equation*}
\mathrm{Y}(Q)=c_{1}\mathrm{M}(Q)+c_{2}\rho_{\pi/ 2}
\end{equation*}
for all $Q\in \mathcal{P}^2$. Here, $\rho_{\pi/2}$ denotes the counter-clockwise rotation  by $\pi/2$ in $\mathbb{R}^2$.
\end{lem}

\begin{lem}\label{lem9}\cite{Wang2}
Suppose that  $n\geq3$. An operator $\mathrm{Y}:\mathcal{P}^n\rightarrow\langle\mathbb{M}^n,+\rangle$ is a simple, 
 continuous, and $\mathrm{SL}(n)$ covariant valuation if and only if there is a real constant $c$ such that
\begin{equation*}
\mathrm{Y}(Q)=c\mathrm{M}(Q)
\end{equation*}
for all $Q\in\mathcal{P}^n$.
\end{lem}


The following lemma is similar to Tsang's result (\cite{Tsang1}, Lemma 3.2), so the proof is omitted.
\begin{lem}\label{lem12}
Let $n\geq2$. 
If  $\Psi:L^{\phi}(\mathbb{R}^n, |x|^2dx)\rightarrow \langle \mathbb{M}^n, + \rangle$ is a valuation, then for all $\alpha_{1},\ldots,\alpha_{l}\in\mathbb{R},$ $l\in\mathbb{N},$ and  pairwise disjoint sets $M_{1},\ldots,M_{l}\in\mathscr{B}$ with finite measure, 
\begin{equation*}
\Psi\left(\sum\limits_{i=1}\limits^{l}\alpha_{i}\chi_{M_{i}}\right)
=\sum\limits_{i=1}\limits^{l}\Psi(\alpha_{i}\chi_{M_{i}}).
\end{equation*} 
\end{lem}

The following lemma is due to Ludwig (\cite{Ludwig6}, Lemma 3), which describes a significant connection between functions on $\mathcal{P}^n$ and on $L^{\phi}(\mathbb{R}^n, |x|^2dx)$.

\begin{lem}\label{lem01}\cite{Ludwig6}
For $Q\in\mathcal{P}^n$ and $\alpha\in\mathbb{R}$, we have $\mathrm{K}(\alpha\chi_Q)=\alpha\mathrm{M}(Q)$.
\end{lem}
Define  $\Psi'(h)=\Psi(h)-\Psi(0)$, then the following lemma holds.

\begin{lem}\label{lem10}
Let $n\geq2$. Suppose that $\phi$ is a Young function satisfying the $\Delta_{2}$-condition. If an operator $\Psi': L^{\phi}(\mathbb{R}^n, |x|^2dx)\rightarrow\langle\mathbb{M}^n, +\rangle$ is a continuous and $\mathrm{SL}(n)$ covariant valuation, then there exists $\xi\in C_{\phi}(\mathbb{R})$  and
\begin{equation*}
\Psi'(h)=\mathrm{K}(\xi\circ h)
\end{equation*}
for any $h\in L^{\phi}_{a}(\mathbb{R}^n, |x|^2dx).$
\end{lem}
\emph{Proof.} 
We divide the proof into three steps as follows. First, we apply Lemma \ref{lem8} and Lemma \ref{lem9} to construct the functional $\xi$. Second, we establish the continuity of $\xi$. 
Finally, we present some properties of $\xi$.

\emph{Construction of $\xi$. } For $\alpha\in\mathbb{R}$, define $\mathrm{Y}_{\alpha}:\mathcal{P}^n\rightarrow\langle\mathbb{M}^n,+\rangle$ by setting
\begin{equation*}
\mathrm{Y}_{\alpha}(Q)=\Psi'(\alpha\chi_{Q})
\end{equation*}
for $Q\in \mathcal{P}^n.$
It is shown that $\mathrm{Y}_{\alpha}$ is an $\mathrm{SL}(n)$ covariant valuation in \cite{Wang2}.
Next, we prove that $\mathrm{Y}_{\alpha}$ is continuous. Applying Lemma \ref{lem2}, there exists a sequence $\{Q_{i}\}$, which is the union of a finite number of closed cubes, such that $\lambda(Q_{i}\bigtriangleup Q)<{1}/{i}$ for $Q\in\mathcal{P}^n$. If $x\in Q_{i}\bigtriangleup Q$, then there is a positive constant $d$ such that $|x|\leq d.$  Observe that
\begin{equation*}
\mu_{n}(Q_{i}\bigtriangleup Q)=\int_{Q_{i}\bigtriangleup Q}|x|^2 dx\leq d^{2}\int_{Q_{i}\bigtriangleup Q}dx= d^{2}\lambda(Q_{i}\bigtriangleup Q)<d^{2}/i,
\end{equation*}
which implies that
\begin{equation*}
  \lim\limits_{i\rightarrow\infty }\mu_{n}(Q_{i}\bigtriangleup Q)=0.  
\end{equation*} If $\mu_{n}(Q_{i}\bigtriangleup Q)=0$, then $\|\chi_{Q_{i}\bigtriangleup Q}\|_{\phi}=0$. Hence, $\|\alpha\chi_{Q_{i}}-\alpha\chi_{Q}\|_{\phi}=0$.  In the other case, if $0<\mu_{n}(Q_{i}\bigtriangleup Q)<\infty$, then
\begin{equation*}
\|\alpha\chi_{Q_{i}}-\alpha\chi_{Q}\|_{\phi}
=|\alpha| \|\chi_{Q_{i}\bigtriangleup Q}\|_{\phi}
=|\alpha| \phi^{*^{-1}}\left(\frac{1}{\mu_{n}(Q_{i}\bigtriangleup Q)}\right)\mu_{n}(Q_{i}\bigtriangleup Q).
\end{equation*}
 Let $t=1\big{/}\mu_{n}(Q_{i}\bigtriangleup Q)$, we obtain
 \begin{equation*}
\|\alpha\chi_{Q_{i}}-\alpha\chi_{Q}\|_{\phi}=|\alpha|\frac{\phi^{*^{-1}}(t)}{t}.
\end{equation*}
 Applying Lemma \ref{lem1}, we have $\|\alpha\chi_{Q_{i}}-\alpha\chi_{Q}\|_{\phi}\rightarrow0$ as $i\rightarrow\infty$. 
Thus, due to the continuity of $\Psi'$,  it leads to that
 $\mathrm{Y}_{\alpha}$ is continuous on $\mathcal{P}^n$.
For $Q\in\mathcal{P}^n$ with $\lambda(Q)=0$, it follows from the $\|\alpha\chi_{Q}\|_{\phi}=0$ that $\alpha\chi_{Q}=0$ a.e. in $L^{\phi}(\mathbb{R}^n, |x|^2dx)$. Thus, $\mathrm{Y}_{\alpha}(Q)=\Psi'(0)=0$, indicating  that $\mathrm{Y}_{\alpha}$ is simple.
By  Lemma \ref{lem9}, there is a real constant  $c_{\alpha}$ such that
\begin{equation}\label{59}
\Psi'(\alpha\chi_{Q})=c_{\alpha}\mathrm{M}(Q)
\end{equation}
for $Q\in\mathcal{P}^n$.
Define the functional $\xi:\mathbb{R}\rightarrow\mathbb{R}$ by setting
\begin{equation*}
\xi(\alpha)=c_{\alpha}
\end{equation*}
for $\alpha\in\mathbb{R}$. 

\emph{Continuity of $\xi$.}
Let $\{\alpha_{i}\}\subset \mathbb{R}$ be a sequence satisfying $\alpha_{i}\rightarrow\alpha$ as $i\rightarrow\infty$. Assume that $\mu_{n}(Q)\neq0$.
Observe that
\begin{equation*}
\|\alpha_{i}\chi_{Q}-\alpha\chi_{Q}\|_{\phi}
=|\alpha_{i}-\alpha|\phi^{*^{-1}}\left(\frac{1}{\mu_{n}(Q)}\right)\mu_{n}(Q).
 \end{equation*}
It follows that $\|\alpha_{i}\chi_{Q}-\alpha\chi_{Q}\|_{\phi}\rightarrow0$
as $\alpha_{i}\rightarrow\alpha$. By the continuity of $\Psi'$, we have $\Psi'(\alpha_{i}\chi_{Q})\rightarrow\Psi'(\alpha\chi_{Q})$. Since 
\begin{equation*}
\xi(\alpha)=\frac{\Psi'(\alpha\chi_{Q})}{\mathrm{M}(Q)} \,\, \mathrm{for}\,\, Q\in\mathcal{P}^n, \, \, \mathrm{M}(Q)\neq0,
\end{equation*} 
 we have  $\xi(\alpha_{i})\rightarrow\xi(\alpha)$ as $\alpha_{i}\rightarrow\alpha$.
Similar to the proof of Lemma 8 in \cite{Tsang2}, we obtain
$\Psi'(s) = \mathrm{K}(\xi \circ s)$, where $s = \sum\limits_{i=1}^{l} \alpha \chi_{M_{i}}$,  $\alpha \in \mathbb{R}$, $M_{i} \subset A[r,R)$ and $s\in L^{\phi}_{a}(\mathbb{R}^n, |x|^2dx)$. The detailed proof will be omitted.

 
 \emph{Properties of $\xi$.} Assume that $\alpha=0$ and $\mathrm{M}(Q)\neq 0$ for $Q\in\mathcal{P}^n$ in (\ref{59}).
 By Lemma \ref{lem4}, we have $\xi(0)=0$.
  For any $\gamma\geq0$, there exists some $\beta_{0}\neq0$ such that $|\xi(\beta_{0})|>\gamma\phi(|\beta_{0}|)$. Then, there is a  sequence $\{\beta_{i}\}$ in $\mathbb{R}\backslash\{0\}$ satisfying
\begin{equation*}
|\xi(\beta_{i})|>2^{i}\phi(|\beta_{i}|).
\end{equation*}
Thus, there exists a subsequence $\{\beta_{i_{j}}\}$ of $\{\beta_{i}\}$ satisfying either
$$\xi(\beta_{i_{j}})>2^{i_{j}}\phi(|\beta_{i_{j}}|)
\quad \text{or} \quad   \xi(\beta_{i_{j}})<-2^{i_{j}}\phi(| \beta_{i_{j}}|).$$
 Let $B_{n}$ denote the unit ball in $\mathbb{R}^n$.
 Tsang (\cite{Tsang2}) constructed a collection of pairwise disjoint balls defined by $\{M_{j}:M_{j}=r_{j}B_{n}+c_{j}e_{1}\}$ in $\mathbb{R}^n$, which satisfy the following conditions: $(1)\ c_{j}>r_{j}$ for $j\geq1$; $(2)\ 0<r_{j}<1$ for $j\geq1$; $(3)\ $ $c_{j}>c_{j-1}+2 $ for $j\geq2$;
 $(4)\ $$(c_{j}+r_{j})^2\lambda(M_{j})=\beta'_{j}$ for $j\geq1$, where $\{\beta_{j}'\}$ is a sequence of positive numbers in $ \mathbb{R}$;$\ $(5)\ $\lim\limits_{j\rightarrow\infty}\frac{c_{j}^2}{(c_{j}+r_{j})^2}=1$. Let $\beta'_{j}=\frac{1}{2^{i_{j}}\phi(|\beta_{i_{j}}|)}$ for all $j\in\mathbb{N}$.
We define $h:\mathbb{R}^n\rightarrow\mathbb{R}$ by
\begin{equation*}
h=\sum\limits_{j=1}\limits^{\infty}\beta_{i_{j}}\chi_{M_{j}}.
\end{equation*}
Since $\phi$ satisfies the $\Delta_{2}$-condition, in order to prove that $h\in L^{\phi}(\mathbb{R}^n, |x|^2dx)$, it is sufficient to show that $\rho_{\phi}(h)<\infty$. In fact, 
\begin{equation*}
\begin{aligned}
\rho_{\phi}(h)
&=\displaystyle\int_{\mathbb{R}^n}\phi\left(\left|\sum\limits_{j=1}\limits^{\infty}
\beta_{i_{j}}\chi_{M_{j}}\right|\right)|x|^2dx\\
&=\sum\limits_{j=1}\limits^{\infty}\phi(|\beta_{i_{j}}|)
\int_{M_{j}}|x|^2 dx\\
&\leq\sum\limits_{j=1}\limits^{\infty}\phi(| \beta_{i_{j}}|)(c_{j}+r_{j})^2\lambda(M_{j})\\
&=\sum\limits_{j=1}\limits^{\infty}\frac{1}{2^{i_{j}}}=1<\infty.
\end{aligned}
\end{equation*}
If $\xi(\beta_{i_{j}})>2^{i_{j}}\phi(|\beta_{i_{j}}|)$ for $j\in\mathbb{N}$, then 
\begin{align*}
\sum\limits_{j=1}\limits^{\infty}\xi(\beta_{i_{j}})c_{j}^2\lambda(M_{j})
&>\sum\limits_{j=1}\limits^{\infty}2^{i_{j}}\phi(| \beta_{i_{j}}|)c_{j}^2\lambda(M_{j})\\
&=\sum\limits_{j=1}\limits^{\infty}\frac{c_{j}^2}{(c_{j}+r_{j})^2}=\lim\limits_{l\rightarrow\infty}\sum\limits_{j=1}\limits^{l}\frac{c_{j}^2}{(c_{j}+r_{j})^2}=\infty.
\end{align*}
Observe that $\int_{M_{j}}x_{i}x_{j}dx=\lambda(M_{j})(c_{j}^2,0,\ldots,0)^{t}$ for $j\in\mathbb{N}$, so 
\begin{align*}
\mathrm{K}_{ij}(\xi\circ h)=\int_{\mathbb{R}^n}(\xi\circ h)(x)x_{i}x_{j}dx&=\left(\sum\limits_{j=1}\limits^{\infty}\xi(\beta_{i_{j}})c_{j}^2\lambda(M_{j}),0,\ldots,0\right)^{t}.
\end{align*}
Hence,
$\left|\mathrm{K}(\xi\circ h)\right|=\infty.$

If  $\xi(\beta_{i_{j}})<-2^{i_{j}}\phi(|\beta_{i_{j}}|)$ for $j\in\mathbb{N}$, then  
\begin{align*}
\sum\limits_{j=1}\limits^{\infty}\xi(\beta_{i_{j}})c_{j}^2\lambda(M_{j})
&<-\infty.
\end{align*} 
A similar argument as in the first case yields
$\left|\mathrm{K}(\xi\circ h)\right|=\infty,$
which contradicts the fact that
\begin{equation*}
\left|\mathrm{K}(\xi\circ h)\right|<\infty.
\end{equation*}
It follows that $|\xi(\beta)|\leq\gamma\phi(|\beta|)$ for $\beta\in\mathbb{R}$. Since simple functions are dense in $L^{\phi}_{a}(\mathbb{R}^n, |x|^2dx)$, then for any $h\in L^{\phi}_{a}(\mathbb{R}^n, |x|^2dx),$ there exists a sequence $\{s_{i}\}$ such that $\|s_{i}-h\|_{\phi} \rightarrow 0$. According to Lemma \ref{lem6}, consider a subsequence $\{s_{k}\}$ of $\{s_{i}\}$ that converges to $ h(x)$ a.e. Due to the continuity of $\xi$, we have
$\lim_{k\rightarrow\infty}(\xi\circ s_{k})(x)x_{i}x_{j}=(\xi\circ h)(x)x_{i}x_{j}, \,\,a. e.$
Since $s_{k}\in L^{\phi}_{a}(\mathbb{R}^n, |x|^2dx)$, it follows that $|\mathrm{K}(\xi\circ  s_{k})|<\infty$. From the above discussion, we have
\begin{equation*}
|(\xi\circ s_{k})(x)x_{i}x_{j}|\leq|(\xi\circ s_{k})(x)||x|^2\leq\gamma\phi(\vert s_{k}(x)\vert)|x|^2.
\end{equation*}
As $\gamma\phi(|s_{k}(x)|)|x|^2$ is integrable, the dominated convergence theorem implies that
\begin{equation*}
\lim\limits_{k\rightarrow\infty}\int_{\mathbb{R}^n}(\xi\circ s_{k})(x)x_{i}x_{j}dx=\int_{\mathbb{R}^n}(\xi\circ h)(x)x_{i}x_{j}dx.
\end{equation*}
Additionally, due to $\|s_{k}-h\|_{\phi} \rightarrow 0$ and the continuity of $\Psi'$, we obtain $\Psi'(s_{k})\rightarrow\Psi'(h)$ as $k\rightarrow\infty$.
This yields
\begin{equation*}
\lim\limits_{k\rightarrow\infty}\Psi'(s_{k})=\lim\limits_{k\rightarrow\infty}\int_{\mathbb{R}^n}(\xi\circ s_{k})(x)x_{i}x_{j}dx=\int_{\mathbb{R}^n}(\xi\circ h)(x)x_{i}x_{j}dx=\Psi'(h)
\end{equation*}
for any $h\in L^{\phi}_{a}(\mathbb{R}^n, |x|^2dx)$.
\qed
\begin{lem}\label{lem01} Let $n\geq2$.
An operator $\Psi':L^{\phi}(\mathbb{R}^n, |x|^2dx)\rightarrow\langle\mathbb{M}^n,+\rangle$ is a valuation and $\xi\in C_{0}(\mathbb{R})$. If  for any non-negative and any non-positive $h\in L^{\phi}(\mathbb{R}^n, |x|^2dx)$, $\Psi'(h)=\mathrm{K}(\xi\circ h)$, then
\begin{equation*}
\Psi'(h)=\mathrm{K}(\xi\circ h)
\end{equation*}
for any $h\in L^{\phi}(\mathbb{R}^n, |x|^2dx)$.
\end{lem}
\emph{Proof.} For $h\in L^{\phi}(\mathbb{R}^n, |x|^2dx),$  note that $h=h\vee0+h\wedge0$, where $h\vee0$ is non-negative and $h\wedge0$ is non-positive in $L^{\phi}(\mathbb{R}^n, |x|^2dx)$.
From the definition of $\mathrm{K}$, we have
\begin{align*}
\mathrm{K}_{ij}(\xi\circ (h\vee0))+\mathrm{K}_{ij}(\xi\circ (h\wedge0))&=\int_{\mathbb{R}^n}\xi\circ ((h\vee0)+ (h\wedge0))x_{i}x_{j}dx\\
&=\int_{\mathbb{R}^n}(\xi\circ h)(x)x_{i}x_{j}dx\\
&=\mathrm{K}_{ij}(\xi\circ h)
\end{align*}
for $h\in L^{\phi}(\mathbb{R}^n, |x|^2dx)$. From $\Psi'(0)=0$, it  shows that
$\Psi'(h)=\Psi'(h)+\Psi'(0)=\Psi'(h\vee0)+\Psi'(h\wedge0)=\mathrm{K}(\xi\circ (h\vee0))+\mathrm{K}(\xi\circ (h\wedge0))=\mathrm{K}(\xi\circ h)$
for $h\in L^{\phi}(\mathbb{R}^n, |x|^2dx)$.
\qed
\begin{lem}\label{lem17}
Let $n\geq2$.
Suppose that $\phi$ is a Young function satisfying the $\Delta_{2}$- condition. If an operator $\Psi': L^{\phi}(\mathbb{R}^n, |x|^2dx)\rightarrow\langle\mathbb{M}^n, +\rangle$ is a continuous and $\mathrm{SL}(n)$ covariant valuation 
, then there is a functional $\xi\in C_{\phi}(\mathbb{R}),$ 
and
\begin{equation}\label{001}
\Psi'(h)=\mathrm{K}(\xi\circ h)
\end{equation}
for any $h\in L^{\phi}(\mathbb{R}^n, |x|^2dx).$
\end{lem}
\emph{Proof.}
By Lemma \ref{lem01}, it suffices to show that  (\ref{001}) holds for $ h \geq 0 $ and  $ h \leq 0 $.
Let $h\in L^{\phi}(\mathbb{R}^n, |x|^2dx)$ such that $h\geq 0$. For  $M_{j}=A[1/2^{j+1},1/2^{j})\cup A[j+1,j+2)$ with $j\geq0$  and $k\in\mathbb{N}$, we define $h_{k}:\mathbb{R}^n\rightarrow \mathbb{R}$ as
\begin{equation*}
h_{k}=\sum\limits_{j=0}\limits^{k}\chi_{M_{j}}h.
\end{equation*}
Observe that $0\leq h_{1}\leq \cdots\leq h_{k}\leq h$ and $h_{k}\rightarrow h$, a.e. Since $\xi$ is continuous, then
$\lim_{k\rightarrow\infty}(\xi\circ h_{k})(x)x_{i}x_{j}=(\xi\circ h)(x)x_{i}x_{j}, \,\,\mathrm{a.e.}$
If $h_{k}\in L^{\phi}_{a}(\mathbb{R}^n, |x|^2dx)$, then $|\Psi'(h_{k})|<\infty$. By Lemma \ref{lem10}, we have $|(\xi\circ h_{k})(x)|\leq\gamma\phi(h_{k}(x))$. It leads to
$
|(\xi\circ h_{k})(x)x_{i}x_{j}|\leq\gamma\phi(h_{k}(x))|x|^2\leq\gamma\phi(h(x))|x|^2, \,\,\mathrm{a. e.}
$
Since $\gamma\phi(h(x))|x|^2$ is integrable, by applying the dominated convergence theorem,
\begin{equation*}
\lim\limits_{k\rightarrow\infty}\int_{\mathbb{R}^n}(\xi\circ h_{k})(x)x_{i}x_{j}dx=\int_{\mathbb{R}^n}(\xi\circ h)(x)x_{i}x_{j}dx.
\end{equation*}
 Since $\phi$ is a increasing function  and $|h_{k}-h|\leq h$, then $\phi(|h_{k}-h|)\leq\phi( h)$. Observe that $\phi( h)$ is integrable. Using the dominated convergence theorem again, we deduce that
\begin{equation*}
\lim\limits_{k\rightarrow\infty}\rho_{\phi}(|h_{k}-h|)=\lim\limits_{k\rightarrow\infty}\int_{\mathbb{R}^n}\phi(|h_{k}-h|)|x|^2dx=0,
\end{equation*}
  which implies that $\|h_{k}-h\|_{\phi}\rightarrow 0$ as $k\rightarrow0$. Due to the continuity of $\Psi'$, we obtain $\Psi'(h_{k})\rightarrow\Psi'(h)$ as $k\rightarrow\infty$. Therefore, 
\begin{equation*}
\Psi'(h)=\int_{\mathbb{R}^n}(\xi\circ (h)(x)x_{i}x_{j}dx.
\end{equation*}
Let $h\in L^{\phi}(\mathbb{R}^n, |x|^2dx)$ with $h\leq 0$. Similarly, we  define 
\begin{equation*}
h_{k}=\sum\limits_{j=0}\limits^{k}\chi_{M_{j}}h,
\end{equation*}
where $M_{j}=A[1/2^{j+1},1/2^{j})\cup A[j+1,j+2)$ for $j\geq0$. 
It implies that $h \leq h_{k}\leq \cdots h_{1}\leq 0$ and $h_{k}\rightarrow h$, a.e. Similar discussion as above, we also obtain
\begin{equation*}
\Psi'(h)=\int_{\mathbb{R}^n}(\xi\circ h(x))x_{i}x_{j} dx.
\end{equation*}

\qed

\begin{theo}\label{theo5}
Suppose that $\phi$ is a Young function satisfying the $\Delta_{2}$-condition. If an operator $\Psi: L^{\phi}(\mathbb{R}^2, |x|^2dx)\rightarrow\langle\mathbb{M}^2, +\rangle$ is a continuous and $\mathrm{SL}(2)$ covariant valuation, then there is a functional $\xi\in C_{\phi}(\mathbb{R})$ and real constant $c$ such that
\begin{equation*}
\Psi(h)=\mathrm{K}(\xi\circ h)+c\rho_{\pi/2}
\end{equation*}
for any $h\in L^{\phi}(\mathbb{R}^2, |x|^2dx)$. Here $\rho_{\pi/2}$ denotes the counter-clockwise rotation  by $\pi/2$ in $\mathbb{R}^2$.
\end{theo}
\emph{Proof.} 
By Lemma \ref{lem3}, Lemma \ref{lem17} and $\Psi'(h)=\Psi(h)-\Psi(0)$, we have
\begin{equation*}
\Psi(h)=\Psi'(h)+\Psi(0)=\Psi'(h)+c\rho_{\pi/2}=\mathrm{K}(\xi\circ h)+c\rho_{\pi/2}
\end{equation*}
for $h\in L^{\phi}(\mathbb{R}^2, |x|^2dx)$.
\qed
\begin{theo}\label{theo6}
Let $n\geq3$. Suppose that $\phi$ is a Young function satisfying the $\Delta_{2}$-condition. If an operator $\Psi: L^{\phi}(\mathbb{R}^n, |x|^2dx)\rightarrow \langle\mathbb{M}^n, +\rangle$ is a continuous and $\mathrm{SL}(n)$ covariant valuation, then there is a functional $\xi\in C_{\phi}(\mathbb{R})$ such that
\begin{equation*}
\Psi(h)=\mathrm{K}(\xi\circ h)
\end{equation*}
for any $h\in L^{\phi}(\mathbb{R}^n, |x|^2dx)$. 
\end{theo}
\emph{Proof.} 
By Lemma \ref{lem4}, Lemma \ref{lem17} and $\Psi'(h)=\Psi(h)-\Psi(0)$, we have
\begin{equation*}
\Psi(h)=\Psi'(h)+\Psi(0)=\Psi'(h)=\mathrm{K}(\xi\circ h)
\end{equation*}
for $h\in L^{\phi}(\mathbb{R}^n, |x|^2dx)$.
\qed

\par \textbf{Acknowledgement}\quad The authors would like to express their gratitude to Monika Ludwig for her advice and assistance with this paper.




\end{document}